\documentclass{amsart}
\usepackage{graphicx}
\usepackage{epstopdf}
\usepackage[colorlinks,allcolors=blue]{hyperref}
\usepackage{amsmath,amssymb,mathtools}

\newtheorem{theorem}{Theorem}[section]
\newtheorem{lemma}{Lemma}[section]

\theoremstyle{remark}
\newtheorem{remark}{Remark}[section]
\newtheorem{notation}{Notation}[section]

\numberwithin{equation}{section}

\DeclareMathOperator{\sinc}{sinc}
\DeclareMathOperator{\sign}{sign}

\usepackage[capitalize,nameinlink]{cleveref}

\Crefname{figure}{Figure}{Figures}

\crefformat{equation}{\textup{#2(#1)#3}}
\crefrangeformat{equation}{\textup{#3(#1)#4--#5(#2)#6}}
\crefmultiformat{equation}{\textup{#2(#1)#3}}{ and \textup{#2(#1)#3}}
{, \textup{#2(#1)#3}}{, and \textup{#2(#1)#3}}
\crefrangemultiformat{equation}{\textup{#3(#1)#4--#5(#2)#6}}%
{ and \textup{#3(#1)#4--#5(#2)#6}}{, \textup{#3(#1)#4--#5(#2)#6}}{, and \textup{#3(#1)#4--#5(#2)#6}}

\crefdefaultlabelformat{#2\textup{#1}#3}

\begin{document}

\title{On the inverse of \lowercase{$\sinc x$}}

\author{Aleksey V. Kargovsky}
\email{kargovsky@yumr.phys.msu.ru}
\address{Faculty of Physics, Lomonosov Moscow State University, Leninskie Gory, 119992 Moscow, Russia}

\begin{abstract}
Inversion of function $\sinc x$ is studied. New series and integral representations of branches of inverse function are obtained using Fourier analysis.
\end{abstract}

\maketitle

\section{Introduction and preliminaries}
The $\sinc$ function, also called \emph{cardinal sine} or \emph{sampling function}, arises frequently in mathematics, physics and engineering in normalized and unnormalized forms. 

Cardinal sine is closely related to the spherical Bessel functions of the first kind $j_n(z)$ \cite{AS} and, in particular,
\begin{equation}
\begin{aligned}
\sinc z &=j_0(z),\\
\sinc^\prime z &=-j_1(z).
\end{aligned}
\end{equation} 

The extrema of $\sinc x$ in real axis are located at the points $\pm j_{3/2,n}$ and zero, where $j_{3/2,n}$ is the $n$th positive root of the Bessel function of the first kind $J_{3/2}(x)$, and the zeros of $\sinc x$ are at non-zero integer multiples of $\pi$.

Here we consider the transcendental equation
\begin{equation}
\frac{\sin y}{y} = x,\label{sinc}
\end{equation}
which apparently first studied in detail by Hardy~\cite{Har02}. He proved that the roots of \cref{sinc}, where $x$ is real and positive, approach asymptotically to the points
\begin{equation*}
\pm\left(2p+\frac{1}{2}\right)\pi \pm i \ln\boldsymbol{(}(4p+1)\pi x\boldsymbol{)},\;p\in\mathbb{N}.
\end{equation*}

On the theoretical basis laid in Hardy's manuscript, further studies were conducted on the numerical analysis of the roots of transcendental equation~\cref{sinc}, its generalizations and analogues \cite{Bur73,Fet76,Han97,Mis64}.

It should be noted that similar equations are encountered in mathematical physics when considering some classes of biharmonic problems in the theory of Stokes flow \cite{Jos77,Kat00,Yoo78} and elasticity \cite{Buc64,Buc65,Fad40,Gay64,Nar65}. 
They also arise in quantum mechanics when considering a wide range of problems from illustrative examples \cite{Mes99} to the description of quantum dots and quantum wires \cite{Bar15}.

Closed-form solution of \cref{sinc} expressed in terms of elementary quadratures was reported in \cite{Bur73} for arbitrary argument. The method used there is based on the solution of a homogeneous Riemann-Hilbert boundary value problem in the theory of
sectionally analytic functions \cite{Bur73b}. However, in some cases, using the resulting expressions may not be very convenient.

In this paper we focus on the study of real solutions of \cref{sinc} and derive series and integral representations of the inverse $\sinc$ function using Fourier analysis.

\section{Main Results}

Using the approach proposed in \cite{Gai23}, consider equation
\begin{equation}
y = \frac{\sin \xi}{x},\;x\in[-1,\,1].
\end{equation}
Using \cref{sinc}, we derive
\begin{equation}
\sin \xi + (-1)^{n+1}x\xi = n\pi x,\;n\in\mathbb{Z}.
\label{sinw}
\end{equation}

Let us consider functions
\begin{equation}
u_\pm(\xi) = \sin \xi \pm x \xi,\;|\xi|\leqslant \varphi_\pm=\arccos (\mp x).
\end{equation}

Since they are continuous and strictly monotonic in given interval, there are strictly monotonic continuous inverse functions \cite{Jef99}:
\begin{equation}
v_\pm(\xi) = u_\pm^{-1}(\xi),\;|\xi|\leqslant l_\pm=u_\pm\boldsymbol{(}\arccos (\mp x)\boldsymbol{)}.\label{vpm}
\end{equation}

Taking into account \cref{sinw}, for $k$th branch of inverse $\sinc$ we can write
\begin{equation}
\sinc_k^{-1} x = \frac{\sin v_\pm(\pi k x)}{x},\;k\in\mathbb{Z}\,\backslash\,\{0\},\label{arcsinc}
\end{equation}
where the upper sign applies to odd branches, and the lower to even ones.

For brevity, we introduce some simplified notations used throughout this paper.
\begin{notation}
For $k\in\mathbb{Z}\,\backslash\,\{0\}$, we make the following definitions:
\begin{enumerate}
\item $j^{(1)}_k \coloneqq j_{\frac{3}{2},\,2\left\lceil\frac{|k|}{2}\right\rceil-1}$ and $j^{(2)}_k \coloneqq j_{\frac{3}{2},\,2\left\lfloor\frac{|k|}{2}\right\rfloor}$;
\item $x^{(1,2)}_k \coloneqq \sinc j^{(1,2)}_k$;
\item $T^{(1,2)}_k \coloneqq \arccos\left((-1)^k x^{(1,2)}_k\right)$;
\end{enumerate}
where $j_{3/2,n}$ is the $n$th positive root of the Bessel function of the first kind $J_{3/2}(x)$ and $j_{3/2,0}=0$.
\end{notation}
Note that the domain of $\sinc^{-1}_k x$ is $[x^{(1)}_k,\,x^{(2)}_k]$.

We first prove the following lemma.
\begin{lemma}
\label{lem1}
Let $k$ is a non-zero integer and $x\in[x^{(1)}_k,\, x^{(2)}_k]$.
Then $l_\pm(x)\geqslant|\pi k x|$, and equality holds only at the end-points of given interval. 
\end{lemma}
\begin{proof}
Consider continuous function $\Delta(x)=l_\pm(x)-|\pi k x|$ at given interval. 
Since
\begin{equation}
\Delta^\prime(x)=\pm\arccos(\mp x)-\pi |k| \sign x,\;x\neq 0,
\end{equation}
$\Delta(x)$ strictly increases for $x<0$ and strictly decreases for $x>0$, and it reaches a maximum $\Delta(0)=1$ at $x=0$.

At points $j^{(1,2)}_k\sign k$ cardinal sine has local extrema and, consequently, we have
\begin{equation}
x^{(1,2)}_k = \frac{\sin j^{(1,2)}_k}{j^{(1,2)}_k}=\cos j^{(1,2)}_k.\label{sincext}
\end{equation}
Namely, there are local minima at points $j^{(1)}_k\sign k$ ($\cos j^{(1)}_k < 0$) and local maxima at points $j^{(2)}_k\sign k$ ($\cos j^{(2)}_k > 0$).

Let us show that 
\begin{equation}
\arccos\left(\mp x^{(1,2)}_k\right)=\mp\left(j^{(1,2)}_k-\pi |k|\right)\sign\left(\cos j^{(1,2)}_k\right).\label{arccoscos}
\end{equation}

\noindent\emph{Case} 1. Let $k$ is odd. Then we have, except $j^{(2)}_1=j^{(2)}_{-1}=0$,
\begin{equation*}
(|k|-1)\pi < j^{(2)}_k < |k|\pi,\quad|k|\pi < j^{(1)}_k < (|k|+1)\pi,
\end{equation*}
or, more precisely \cite{Wat95},
\begin{equation*}
\left(|k|-\frac{3}{4}\right)\pi < j^{(2)}_k < \left(|k|-\frac{1}{2}\right)\pi,\quad
\left(|k|+\frac{1}{4}\right)\pi < j^{(1)}_k < \left(|k|+\frac{1}{2}\right)\pi.
\end{equation*}
Therefore we get
\begin{align}
\nonumber&\arccos\left(-x^{(2)}_k\right)=\arccos\left(-\cos j^{(2)}_k\right)>\frac{\pi}{2} \Rightarrow \arccos\left(-\cos j^{(2)}_k\right)=|k|\pi-j^{(2)}_k,\\
\nonumber&\arccos\left(-x^{(1)}_k\right)=\arccos\left(-\cos j^{(1)}_k\right)<\frac{\pi}{2} \Rightarrow \arccos\left(-\cos j^{(1)}_k\right)=j^{(1)}_k-|k|\pi,\\
&\arccos\left(-x^{(1,2)}_k\right)=-\left(j^{(1,2)}_k-|k|\pi\right)\sign\left(\cos j^{(1,2)}_k\right).\label{arccoscos1}
\end{align}
In the special case $j^{(2)}_1=j^{(2)}_{-1}=0$ and $x^{(2)}_1=x^{(2)}_{-1}=1$, we obtain
\begin{equation}
\arccos\left(-x^{(2)}_1\right)=\arccos\left(-x^{(2)}_{-1}\right)=\pi,
\end{equation}
and \cref{arccoscos1} holds.

\noindent\emph{Case} 2. Let $k$ is even. Then we have
\begin{equation*}
\left(|k|-\frac{3}{4}\right)\pi < j^{(1)}_k < \left(|k|-\frac{1}{2}\right)\pi,\quad
\left(|k|+\frac{1}{4}\right)\pi < j^{(2)}_k < \left(|k|+\frac{1}{2}\right)\pi.
\end{equation*}
Therefore we get
\begin{align}
\nonumber&\arccos\left(x^{(1)}_k\right)=\arccos\left(\cos j^{(1)}_k\right)>\frac{\pi}{2} \Rightarrow \arccos\left(\cos j^{(1)}_k\right)=|k|\pi-j^{(1)}_k,\\
\nonumber&\arccos\left(x^{(2)}_k\right)=\arccos\left(\cos j^{(2)}_k\right)<\frac{\pi}{2} \Rightarrow \arccos\left(\cos j^{(2)}_k\right)=j^{(2)}_k-|k|\pi,\\
&\arccos\left(x^{(1,2)}_k\right)=\left(j^{(1,2)}_k-|k|\pi\right)\sign\left(\cos j^{(1,2)}_k\right).\label{arccoscos2}
\end{align}
From \cref{arccoscos1,arccoscos2} the statement \cref{arccoscos} follows.

Using \cref{sincext,arccoscos}, we obtain
\begin{equation}
\Delta\left(x^{(1,2)}_k\right)=\left|\sin j^{(1,2)}_k\right|-\left|\cos j^{(1,2)}_k\right|\left(j^{(1,2)}_k-\pi |k|\right)-\pi \left|k \cos j^{(1,2)}_k\right| = 0.
\end{equation}
Therefore we have $\Delta(x)\geqslant 0$ as $x\in[x^{(1)}_k,\,x^{(2)}_k]$, where equality holds only at the end-points of given interval. Hence the result. 
\end{proof}

Series representation of inverse $\sinc$ is given in the following theorem.
\begin{theorem}
\label{thm1}
Let $k$ is a non-zero integer and $x\in[x^{(1)}_k,\, x^{(2)}_k]$.
Then $k$th branch of inverse $\sinc$ can be represented as
\begin{equation}
\label{arcsinc1}%
\sinc_k^{-1} x = \frac{\pi k}{l}\sqrt{1-x^2}
+ (-1)^k \sum\limits_{n=1}^\infty\frac{2}{n}  
 A_{\frac{\pi n x}{l}}\boldsymbol{\Bigl(}\arccos\boldsymbol{(}(-1)^k x\boldsymbol{)},\;(-1)^k\frac{\pi n}{l}\boldsymbol{\Bigr)}\sin \frac{\pi^2 n k x}{l},
\end{equation}
where $A_\nu(\varphi,x)$ is the incomplete Anger function and
\begin{equation*}
l=\sqrt{1-x^2}-(-1)^k x \arccos\boldsymbol{(}(-1)^k x\boldsymbol{)}.
\end{equation*}
\end{theorem}
\begin{proof}
Since $u_\pm(\xi)$ is an odd function, its inverse $v_\pm(\xi)$ is also odd. Let us extend function $\sin v_\pm(\xi)$ defined at $|\xi|\leqslant l_\pm$ periodically over $\mathbb{R}$ (with period $2l_\pm$). This extension have discontinuities of the first kind at points $\{(2n+1)l_\pm: n\in\mathbb{Z}\}$, except the case $x=1$. Since $\sin v_\pm(\xi)$ is of bounded variation over $[-l_\pm,l_\pm]$, its Fourier series converges to
\begin{equation}
\frac{1}{2}\left[\sin v_\pm(\xi+0)+\sin v_\pm(\xi-0)\right] = \sum\limits_{n=1}^\infty{b_n \sin \frac{\pi n \xi}{l_\pm}},\label{sinvs}
\end{equation}
moreover, the Fourier series converges uniformly to $\sin v_\pm(\xi)$ in any interval interior to $(-l_\pm,\,l_\pm)$ \cite{Bar64}.
Fourier coefficients in \cref{sinvs} are equal to
\begin{multline}\label{bn}
b_n = \frac{2}{l_\pm}\int_{0}^{l_\pm}{\sin v_\pm(\xi)\;\sin \frac{\pi n \xi}{l_\pm}\;d\xi} 
=\frac{2}{l_\pm}\int_{0}^{\varphi_\pm}{\sin t\;\sin \frac{\pi n u_\pm(t)}{l_\pm}\;du_\pm(t)} \\
=\frac{2}{\pi n}\left[-\left.\sin t\;\cos\frac{\pi n u_\pm(t)}{l_\pm}\right|_0^{\varphi_\pm}
+\int_{0}^{\varphi_\pm}{\cos t\;\cos \frac{\pi n u_\pm(t)}{l_\pm}\;dt}\right] \\
=\frac{2}{\pi n}\left[(-1)^{n-1}\sin \varphi_\pm
+\int_0^{\varphi_\pm}{\cos t\;\cos\left( \frac{\pi n}{l_\pm}[\sin t\pm x t]\right)\;dt}\right]\\
=\frac{2}{\pi n}\left\{(-1)^{n-1}\sin \varphi_\pm+\frac{\pi}{2}\left[A_{\mp\frac{\pi n x}{l_\pm}+1}\left(\varphi_\pm,\frac{\pi n}{l_\pm}\right)
+A_{\mp\frac{\pi n x}{l_\pm}-1}\left(\varphi_\pm,\frac{\pi n}{l_\pm}\right)\right]\right\},
\end{multline}
where $A_\nu(\varphi,x)$ is the incomplete Anger function \cite{Agr71}:
\begin{equation}
A_\nu(\varphi,x) = \frac{1}{\pi}\int_0^{\varphi}{\cos(x \sin t - \nu t)\;dt}.\label{A}
\end{equation}

Substituting \cref{bn,sinvs} into \cref{arcsinc} and summing using \cite{Gra07}, we get
\begin{multline}
S_k(x) = \frac{1}{x}\sum\limits_{n=1}^\infty{b_n \sin \frac{\pi^2 n k x}{l_\pm}}\\
=\frac{\pi k}{l_\pm}\sin\varphi_\pm \left(1-\delta_{l_\pm,|k \pi x|}\right)
\mp \sum\limits_{n=1}^\infty{\frac{2}{n}  A_{\mp\frac{\pi n x}{l_\pm}}\left(\varphi_\pm,\frac{\pi n}{l_\pm}\right)\sin \frac{\pi^2 n k x}{l_\pm}}\\
=\frac{\pi k}{l_\pm}\sqrt{1-x^2} \left(1-\delta_{x,x_k^{(1)}}-\delta_{x,x_k^{(2)}}\right)
\mp \sum\limits_{n=1}^\infty{\frac{2}{n}  A_{\mp\frac{\pi n x}{l_\pm}}\left(\varphi_\pm,\frac{\pi n}{l_\pm}\right)\sin \frac{\pi^2 n k x}{l_\pm}}.\label{sing}
\end{multline}
Here we used \cref{lem1} and recurrence relation for incomplete Anger functions~\cite{Agr71}:
\begin{equation}\label{Arel}
A_{\nu+1}(\varphi,x)+A_{\nu-1}(\varphi,x) = \frac{2\nu}{x}A_{\nu}(\varphi,x) + \frac{2}{\pi x}\sin(x\sin \varphi-\nu\varphi)\qquad(x,\nu\in\mathbb{R}).
\end{equation}

Obviously, we have
\begin{equation}
S_k(x) = \begin{cases}
\sinc_k^{-1} x &\left(x^{(1)}_k <x< x^{(2)}_k\right);\\
0 &\left(x=x^{(1,2)}_k\right).
\end{cases}\label{sk2}
\end{equation}

At the end-points $x_k^{(1,2)}$ owing to \cref{lem1}, we can write
\begin{multline}\label{skb}
\left.\frac{\pi k}{l_\pm}\sqrt{1-x^2}
\mp \sum\limits_{n=1}^\infty\frac{2}{n}  
 A_{\mp\frac{\pi n x}{l_\pm}}\boldsymbol{\Bigl(}\varphi_\pm,\;\frac{\pi n}{l_\pm}\boldsymbol{\Bigr)}\sin \frac{\pi^2 n k x}{l_\pm}\right|_{x=x_k^{(1,2)}}\\
 =\frac{|\sin j_k^{(1,2)}|}{|\cos j_k^{(1,2)}|}\sign k=j_k^{(1,2)}\sign k = \sinc_k^{-1} x_k^{(1,2)}.
\end{multline}

The result now follows from the first case of \cref{sk2} and \cref{skb}.
\end{proof}

\begin{figure}[!htbp]
\centering
\includegraphics[width=0.6\textwidth]{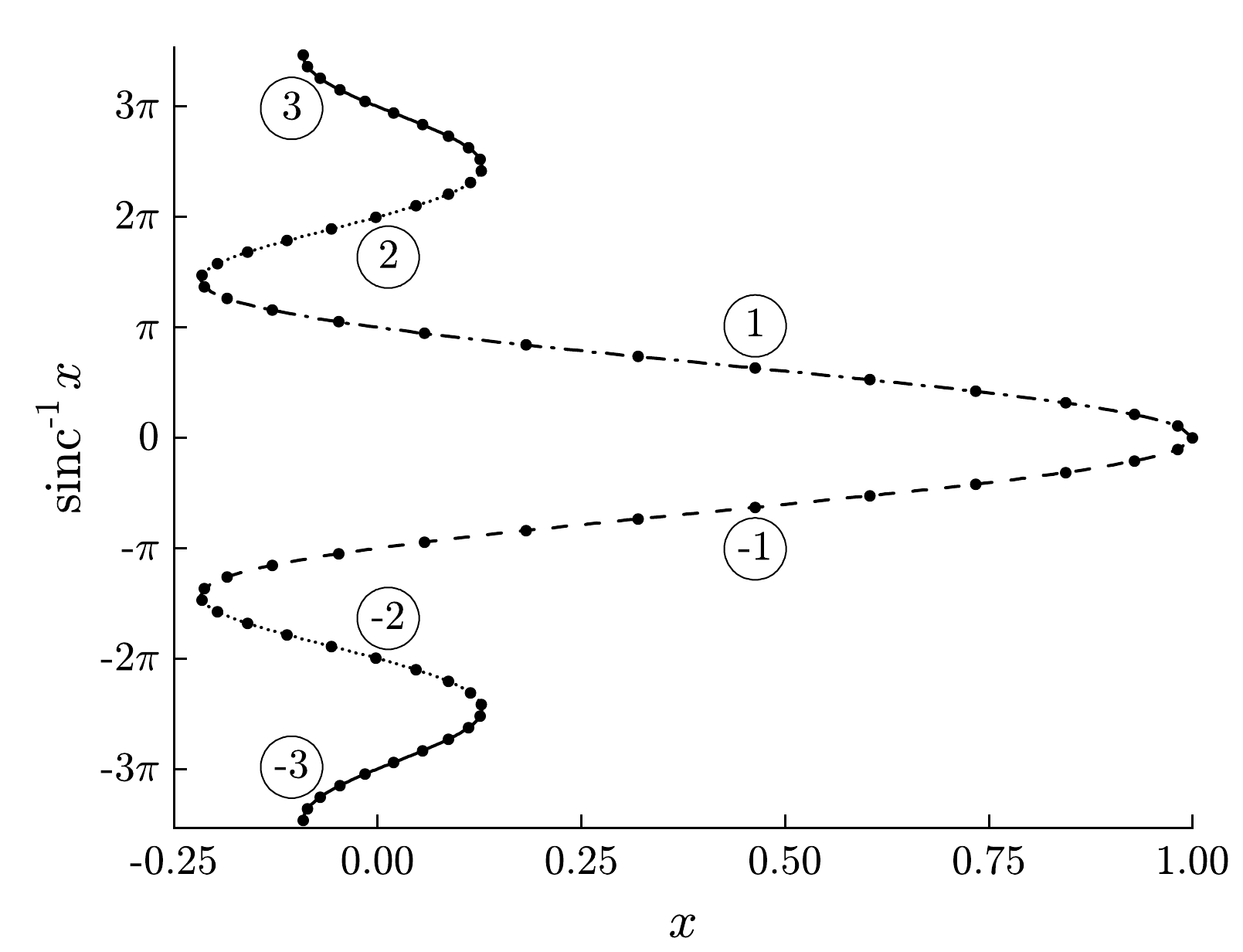}
\caption{Plot of $\sinc^{-1} x$.\label{sinc-1}}
\end{figure}

First three positive and negative branches of function $\sinc^{-1} x$ from \cref{arcsinc1} are presented in \cref{sinc-1}, where the lines represent \cref{arcsinc1} and the markers correspond to inverse mapping $\sinc x \mapsto x$. The numbers in circles denote the branch indexes.

An alternative series representation is given in the following theorem.
\begin{theorem}
\label{thm2}
Let $k$ is a non-zero integer and $x\in[x^{(1)}_k,\, x^{(2)}_k]$.
Then $k$th branch of inverse $\sinc$ can be represented as
\begin{multline}
\label{arcsinc2}%
\sinc_k^{-1} x = \pi k +\sign k \;\Biggl(\frac{|x|}{l}\left[1-\frac{1}{2}\arccos^2\boldsymbol{(}(-1)^k x\boldsymbol{)}\right]\\
+\sign x\Biggl\{\frac{(-1)^k}{l}\left[\sqrt{1-x^2}\arccos\boldsymbol{(}(-1)^k x\boldsymbol{)}-1\right]\\
+ \sum\limits_{n=1}^\infty\frac{2}{n}  
 B_{\frac{\pi n x}{l}}\boldsymbol{\Bigl(}\arccos\boldsymbol{(}(-1)^k x\boldsymbol{)},\;(-1)^k\frac{\pi n}{l}\boldsymbol{\Bigr)}\cos \frac{\pi^2 n k x}{l}\Biggr\}\Biggr),
\end{multline}
where $B_\nu(\varphi,x)$ is the incomplete Weber function and
\begin{equation*}
l=\sqrt{1-x^2}-(-1)^k x \arccos\boldsymbol{(}(-1)^k x\boldsymbol{)}.
\end{equation*}

\end{theorem}
\begin{proof}
First, note that $\sin v_\pm(\xi)=|\sin v_\pm(\xi)|\, \sign \xi$ for $|\xi|\leqslant l_\pm$. Let us extend function $|\sin v_\pm(\xi)|$ defined at $|\xi|\leqslant l_\pm$ periodically over $\mathbb{R}$ (with period $2l_\pm$). This extension is even and continuous, and $|\sin v_\pm(\xi)|$ is of bounded variation over $[-l_\pm,l_\pm]$. Then its Fourier series converges uniformly to it in $\mathbb{R}$ \cite{Bar64}
\begin{equation}
\left|\sin v_\pm(\xi)\right| = \frac{a_0}{2}+\sum\limits_{n=1}^\infty{a_n \cos \frac{\pi n \xi}{l_\pm}}.\label{sinvc}
\end{equation}
Fourier coefficients in \cref{sinvc} are equal to
\begin{multline}\label{a0}
a_0 = \frac{2}{l_\pm}\int_{0}^{l_\pm}{\sin v_\pm(\xi)\;d\xi} 
=\frac{2}{l_\pm}\int_{0}^{\varphi_\pm}{\sin t\;du_\pm(t)} \\
=\frac{2}{l_\pm}\left[\sin t\;u_\pm(t)\biggm|_0^{\varphi_\pm}
-\int_{0}^{\varphi_\pm}{\cos t\;\left(\sin t \pm xt\right)\;dt}\right] \\
=\frac{2}{l_\pm}\left\{l_\pm \sin \varphi_\pm-\left[\frac{1}{2}\sin^2 \varphi_\pm \pm x\left(\varphi_\pm \sin \varphi_\pm +\cos \varphi_\pm -1\right)\right]\right\}\\
=\frac{2}{l_\pm}\left[\frac{1}{2}\sin^2 \varphi_\pm \mp x\left(\cos \varphi_\pm -1\right)\right]
= \frac{2}{l_\pm}\left(\frac{1+x^2}{2}\pm x\right),
\end{multline}
\begin{multline}\label{an}
a_n = \frac{2}{l_\pm}\int_{0}^{l_\pm}{\sin v_\pm(\xi)\;\cos \frac{\pi n \xi}{l_\pm}\;d\xi} 
=\frac{2}{l_\pm}\int_{0}^{\varphi_\pm}{\sin t\;\cos \frac{\pi n u_\pm(t)}{l_\pm}\;du_\pm(t)} \\
=\frac{2}{\pi n}\left[\left.\sin t\;\sin\frac{\pi n u_\pm(t)}{l_\pm}\right|_0^{\varphi_\pm}
-\int_{0}^{\varphi_\pm}{\cos t\;\sin \frac{\pi n u_\pm(t)}{l_\pm}\;dt}\right] \\
=\frac{2}{\pi n}\int_0^{\varphi_\pm}{\cos t\;\sin\left( \frac{\pi n}{l_\pm}[\mp x t-\sin t]\right)\;dt}\\
=\frac{1}{n}\left[B_{\mp\frac{\pi n x}{l_\pm}+1}\left(\varphi_\pm,\frac{\pi n}{l_\pm}\right)
+B_{\mp\frac{\pi n x}{l_\pm}-1}\left(\varphi_\pm,\frac{\pi n}{l_\pm}\right)\right]\\
=\frac{2}{n}\left\{\mp x B_{\mp\frac{\pi n x}{l_\pm}}\left(\varphi_\pm,\frac{\pi n}{l_\pm}\right)
+\frac{l_\pm}{\pi^2 n}\left[(-1)^n-1\right]\right\},\;n\geqslant 1.
\end{multline}
Here $B_\nu(\varphi,x)$ is the incomplete Weber function \cite{Agr71}:
\begin{equation}
B_\nu(\varphi,x) = \frac{1}{\pi}\int_0^{\varphi}{\sin(\nu t - x \sin t)\;dt},\label{B}
\end{equation}
and we used for them recurrence relation suitable for real $x$ and $\nu$:
\begin{equation}\label{Brel}
B_{\nu+1}(\varphi,x)+B_{\nu-1}(\varphi,x) = \frac{2\nu}{x}B_{\nu}(\varphi,x) + \frac{2}{\pi x}\left[\cos(x\sin \varphi-\nu\varphi)-1\right].
\end{equation}

Collecting together \cref{a0,an,sinvc} and summing using \cite{Gra07}, we get
\begin{multline}
\sinc^{-1}_k x = \frac{\left|\sin v_\pm(\pi k x)\right|}{|x|}\sign k
=\pi k +\sign k \;\Biggl\{\frac{|x|}{l_\pm}\left(1-\frac{\varphi_\mp^2}{2}\right)\pm\sign x\\
\times\left[\frac{1}{l_\pm}\left(1-\varphi_\pm\sqrt{1-x^2}\right)
-\sum\limits_{n=1}^\infty\frac{2}{n}  
 B_{\mp\frac{\pi n x}{l}}\left(\varphi_\pm,\;\frac{\pi n}{l_\mp}\right)\cos \frac{\pi^2 n k x}{l_\pm}\right]\Biggr\}.\label{sing2}
\end{multline}

Using \cite{Pru90} it can be verified that the expression in square brackets in \cref{sing2} is zero at $x=0$ and therefore there is no discontinuity.

Hence the result.
\end{proof}

\begin{remark}
Putting $k=1$ and $x=1$ in \cref{arcsinc2}, we can derive the sum of a series containing the Weber functions
\begin{equation}
\sum_{n=1}^\infty{\frac{(-1)^n}{n}\mathbf{E}_{-n}(n)}=\frac{\pi}{4}+\frac{1}{\pi}.
\end{equation}
\end{remark}


Next we obtain integral representation of inverse $\sinc$ using the Fourier's sine formula.

We require the following lemmas.

\begin{lemma}
\label{lem2}
Let $k$ is a non-zero integer and $t\in[0,\,T^{(1,2)}_k]$.
Then 
\begin{equation*}
0\leqslant\sin t \pm t x^{(1,2)}_k\leqslant|\pi k x^{(1,2)}_k|, 
\end{equation*}
and the lower and upper bounds are attained only at the left and right end-points of the interval, respectively. 
\end{lemma}
\begin{proof}
The function $U_\pm(t)=\sin t \pm t x^{(1,2)}_k$ is continuous and strictly increasing in given interval.

Then $U_\pm(0) < U_\pm(t) < U_\pm(T^{(1,2)}_k)$ at the interior points of this interval, and its lower and upper bounds are attained at the left and right end-points, respectively. Using $U_\pm(0)=0$ and $U_\pm(T^{(1,2)}_k)=l_\pm(x^{(1,2)}_k)=|\pi k x^{(1,2)}_k|$ from \cref{lem1}, the result follows.
\end{proof}

\begin{lemma}
\label{lem3}
It holds that
\begin{multline}
\int_0^{T^{(1,2)}_k}{dt\int_0^\infty{\frac{\sin \omega \pi k  x^{(1,2)}_k}{\omega}\cos\left(\omega\left[\sin t \pm tx^{(1,2)}_k\right]\right)\;d\omega}}\\
= \int_0^\infty{\frac{\sin \omega \pi k  x^{(1,2)}_k}{\omega}\;d\omega\int_0^{T^{(1,2)}_k}{\cos\left(\omega\left[\sin t \pm t x^{(1,2)}_k\right]\right)\;dt}}.
\end{multline}
\end{lemma}
\begin{proof}
Following \cite{Whi96}, the integral 
\begin{equation*}
\int_0^\infty{\frac{\sin \alpha x\;\cos \beta x}{x}\;dx}\qquad(\alpha,\beta\in\mathbb{R})
\end{equation*}
converges uniformly in any closed domain of values of $\alpha$ from which the points $\alpha=\pm\beta$ are excluded.

Then taking info account \cref{lem2}, we conclude that the integral 
\begin{equation}
\int_0^\infty{\frac{\sin \omega \pi k  x^{(1,2)}_k}{\omega}\cos\left(\omega\left[\sin t \pm tx^{(1,2)}_k\right]\right)\;d\omega}
\end{equation}
is uniformly convergent over any subinterval of $[0, T^{(1,2)}_k]$ not ending at $T^{(1,2)}_k$. We may therefore invert the integral
\begin{equation}
\int_0^{T^{(1,2)}_k-\epsilon}{dt\int_0^\infty{\frac{\sin \omega \pi k  x^{(1,2)}_k}{\omega}\cos\left(\omega\left[\sin t \pm tx^{(1,2)}_k\right]\right)\;d\omega}}
\end{equation}
for every $0<\epsilon<T^{(1,2)}_k$.

Since
\begin{multline}
\lim_{\epsilon\to 0}\int_{T^{(1,2)}_k-\epsilon}^{T^{(1,2)}_k}{dt\int_0^\infty{\frac{\sin \omega \pi k  x^{(1,2)}_k}{\omega}\cos\left(\omega\left[\sin t \pm tx^{(1,2)}_k\right]\right)\;d\omega}}\\
=\lim_{\epsilon\to 0}\frac{\pi\epsilon}{2}\sign kx^{(1,2)}_k=0, 
\end{multline}
to prove the statement of the lemma it is sufficient to prove that
\begin{subequations}
\begin{equation}\label{lim}
\lim_{\epsilon\to 0}\int_0^\infty{\frac{\sin \omega \pi k  x^{(1,2)}_k}{\omega}\;d\omega\int_{T^{(1,2)}_k-\epsilon}^{T^{(1,2)}_k}{\cos\left(\omega\left[\sin t \pm t x^{(1,2)}_k\right]\right)\;dt}}=0. \tag{\theparentequation}
\end{equation}

Let us split the outer integral in \cref{lim} into a sum of three integrals over the intervals $[0,1/\pi]$, $[1/\pi,1/\epsilon]$ and $[1/\epsilon,\infty)$, and estimate each of them:

\begin{multline}\label{lim1}%
\left|\int_0^{1/\pi}{\frac{\sin \omega \pi k  x^{(1,2)}_k}{\omega}\;d\omega\int_{T^{(1,2)}_k-\epsilon}^{T^{(1,2)}_k}{\cos\left(\omega\left[\sin t \pm t x^{(1,2)}_k\right]\right)\;dt}}\right|\\
\leqslant \epsilon\int_0^{1/\pi}{\left|\sin \omega \pi k  x^{(1,2)}_k\right|\frac{d\omega}{\omega}} = O(\epsilon),
\end{multline}
\begin{multline}\label{lim2}%
\left|\int_{1/\pi}^{1/\epsilon}{\frac{\sin \omega \pi k  x^{(1,2)}_k}{\omega}\;d\omega\int_{T^{(1,2)}_k-\epsilon}^{T^{(1,2)}_k}{\cos\left(\omega\left[\sin t \pm t x^{(1,2)}_k\right]\right)\;dt}}\right|\\
\leqslant \epsilon \int_{1/\pi}^{1/\epsilon}{\frac{d\omega}{\omega}} = \epsilon\ln\frac{\pi}{\epsilon},
\end{multline}
\begin{multline}\label{lim3}%
\left|\int_{1/\epsilon}^\infty{\frac{\sin \omega \pi k  x^{(1,2)}_k}{\omega}\;d\omega\int_{T^{(1,2)}_k-\epsilon}^{T^{(1,2)}_k}{\cos\left(\omega\left[\sin t \pm t x^{(1,2)}_k\right]\right)\;dt}}\right|\\
=\pi\left|\int_{1/\epsilon}^\infty{\left[A_{\mp\omega x^{(1,2)}_k}\left(T^{(1,2)}_k,\omega\right)-A_{\mp\omega x^{(1,2)}_k}\left(T^{(1,2)}_k-\epsilon,\omega\right)\right]\frac{\sin \omega \pi k  x^{(1,2)}_k}{\omega}\;d\omega}\right| \\
\leqslant \int_{1/\epsilon}^\infty{O(\omega^{-4/3})\;d\omega}=O(\epsilon^{1/3}).
\end{multline}
\end{subequations}
In \cref{lim3} we used that at large $\omega$ the incomplete Anger function decreases no slower than $O\left(\omega^{-1/3}\right)$ \cite{Agr71}.

By using \cref{lim1,lim2,lim3} and applying l'H\^{o}pital's rule, limit \cref{lim} is justified, and this ends the proof. 
\end{proof}

\begin{theorem}
\label{thm3}
Let $k$ is a non-zero integer and $x\in[x^{(1)}_k,\; x^{(2)}_k]$.
Then $k$th branch of inverse $\sinc$ can be represented as
\begin{equation}
\label{arcsinci}%
\sinc_k^{-1} x = \pi k
+ 2(-1)^k \int_0^\infty{A_{\omega x}\boldsymbol{\Bigl(}\arccos\boldsymbol{(}(-1)^k x\boldsymbol{)},\;(-1)^k\omega\boldsymbol{\Bigr)} \frac{\sin \omega\pi k x}{\omega}\;d\omega},
\end{equation}
where $A_\nu(\varphi,x)$ is the incomplete Anger function.
\end{theorem}
\begin{proof}

Consider function $v_\pm(\xi)$ defined by \cref{vpm} at $|\xi|\leqslant l_\pm$ and equal to zero elsewhere. Since function $\sin v_\pm(\xi)$ belongs to $L^1(\mathbb{R})$ and is of bounded variation over $\mathbb{R}$, using~\cite{Tit86}, we have
\begin{equation}
\frac{1}{2}\left[\sin v_\pm(\xi+0)+\sin v_\pm(\xi-0)\right] = \frac{1}{\pi}\int_0^\infty{G_\pm(\omega) \sin \omega \xi\;d\omega},\label{sinvf}
\end{equation}
where
\begin{multline}
G_\pm(\omega) = 2\int_{0}^{l_\pm}{\sin v_\pm(\xi)\;\sin \omega \xi\;d\xi} 
=2\int_{0}^{\varphi_\pm}{\sin t\;\sin \omega u_\pm(t)\;du_\pm(t)} \\
=\frac{2}{\omega}\left[-\sin t\;\cos\omega u_\pm(t)\biggl|_0^{\varphi_\pm}
+\int_{0}^{\varphi_\pm}{\cos t\;\cos \omega u_\pm(t)\;dt}\right] \\
=\frac{2}{\omega}\left[-\sin \varphi_\pm \cos\omega l_\pm
+\frac{\pi}{2}\biggl\{A_{\mp \omega x+1}\left(\varphi_\pm,\omega\right)
+A_{\mp \omega x-1}\left(\varphi_\pm,\omega\right)\biggr\}\right] \\
=\frac{2}{\omega}\left[-\sin \varphi_\pm \cos\omega l_\pm+\frac{1}{\omega}\sin\omega l_\pm
\mp\pi x A_{\mp \omega x}\left(\varphi_\pm,\omega\right)\right].\label{G}
\end{multline}
Moreover, the integral \cref{sinvf} converges uniformly to $\sin v_\pm(\xi)$ in any interval interior to $(-l_\pm,\,l_\pm)$ \cite{Tit86}.

Substituting~\cref{G,sinvf} into \cref{arcsinc} and summing using \cite{Gra07}, we write the integral
\begin{multline}
I_k(x) = \frac{1}{\pi}\int_0^\infty G_\pm(\omega) \frac{\sin \omega\pi k x}{x} \;d\omega \\
=\frac{1}{2 x}\left(-\sin \varphi_\pm \sign kx\;\delta_{l_\pm,|\pi k x|}+2\pi k x\right)
\mp 2 \int_0^\infty{A_{\mp\omega x}\left(\varphi_\pm,\omega\right)\frac{\sin \omega\pi k x}{\omega}\;d\omega}.\label{ik}
\end{multline}
Using \cref{lem1}, \cref{ik} becomes
\begin{multline}
I_k(x) =\pi k
\mp 2 \int_0^\infty{A_{\mp\omega x}\left(\varphi_\pm,\omega\right)\frac{\sin \omega\pi k x}{\omega}\;d\omega}\\
-\frac{1}{2}\left( j^{(1)}_k\delta_{x,\;x^{(1)}_k}
+j^{(2)}_k\delta_{x,\;x^{(2)}_k}\right)\sign k.\label{ik2}
\end{multline}

Consider the second term of \cref{ik2} at points $x^{(1,2)}_k$. Integrating using \cite{Gra07} and \cref{lem2} and taking into account \cref{arccoscos},  we obtain
\begin{multline}
\int_0^\infty{A_{\mp\omega x}\left(T^{(1,2)}_k,\,\omega\right)\frac{\sin \omega \pi k  x^{(1,2)}_k}{\omega}\;d\omega}\\
=\frac{1}{\pi}\int_0^\infty{d\omega\;\frac{\sin \omega \pi k  x^{(1,2)}_k}{\omega}\int_0^{T^{(1,2)}_k}{\cos\left(\omega\left[\sin t \pm t x^{(1,2)}_k\right]\right)\;dt}}\\
=\frac{1}{\pi}\int_0^{T^{(1,2)}_k}{dt\int_0^\infty{\frac{\sin \omega \pi k  x^{(1,2)}_k}{\omega}\cos\left(\omega\left[\sin t \pm tx^{(1,2)}_k\right]\right)\;d\omega}}\\
=\frac{1}{2}T^{(1,2)}_k\sign\left(k  x^{(1,2)}_k\right)
=\mp\frac{1}{2}\left(j^{(1,2)}_k\sign k-\pi k\right).\label{Aint}
\end{multline}
The inversion is justified by \cref{lem3}.

Then from \cref{ik2} and \cref{Aint}, we have
\begin{equation}
I_k(x) = \begin{cases}
\sinc_k^{-1} x &\left(x^{(1)}_k <x< x^{(2)}_k\right);\\
\frac{1}{2}j^{(1,2)}_k\sign k=\frac{1}{2}\sinc_k^{-1}x &\left(x=x^{(1,2)}_k\right).
\end{cases}\label{ik3}
\end{equation}
The last case reflects the fact that in \cref{sinvf} we used the Fourier's sine formula for a function that has discontinuities of the first kind at points $-l_\pm$ and $l_\pm$.

Thus at the interior points of given interval the statement \cref{arcsinci} follows from the first case of \cref{ik3}. 

Using \cref{Aint}, we have at the end-points $x_k^{(1,2)}$
\begin{multline*}
\pi k
+ 2(-1)^k \int_0^\infty{A_{\omega x_k^{(1,2)}}\boldsymbol{\Bigl(}T^{(1,2)}_k,\,(-1)^k\omega\boldsymbol{\Bigr)} \frac{\sin \omega \pi k x_k^{(1,2)}}{\omega}\;d\omega}\\
=j^{(1,2)}_k\sign k = \sinc_k^{-1} x_k^{(1,2)},
\end{multline*}
and the result follows.
\end{proof}

\begin{remark}
Putting $k=1$ and $x=1$ in \cref{arcsinci}, we can derive the integral containing the Anger function
\begin{equation}
\int_0^\infty{\mathbf{J}_{\omega}(-\omega) \frac{\sin \pi\omega}{\omega}\;d\omega} = \frac{\pi}{2}.
\end{equation}
\end{remark}

\section{Applications}
In some applications, such as signal processing and spectroscopy, it is required to know the peak width of a $\sinc x$ or $\sinc^2 x$ functions at a certain level. This can be done using \cref{arcsinc1}, \cref{arcsinc2} or \cref{arcsinci}.
For example, it is easy to obtain FWHM of peaks of $\sinc^2 x$ as the difference between values of two adjacent branches of \cref{arcsinc1} at certain point:
\begin{subequations}
\begin{align}
\Delta\omega_0 &= 2\sinc_{1}^{-1}\left(\frac{1}{\sqrt{2}}\right),\\
\Delta\omega_m &= \sinc_{m+1}^{-1}\left(\frac{1}{\sqrt{2}}\sinc j_{3/2,m}\right)
-\sinc_{m}^{-1}\left(\frac{1}{\sqrt{2}}\sinc j_{3/2,m}\right),\;m>0.\label{dwk}
\end{align}
\end{subequations}
Here $m$ is the peak number.
For large $m$ \cref{dwk} becomes
\begin{equation}
\lim_{m\to\infty}\Delta\omega_m = \pi - 2 \lim_{m\to\infty}\arcsin\left(\frac{\pi m}{\sqrt{2}}|\sinc j_{3/2,m}|\right)
=\frac{\pi}{2}.
\end{equation}

\section{Conclusion}
In this paper we obtain real-valued series and integral representations \cref{arcsinc1,arcsinc2,arcsinci} of $\sinc_k^{-1} x$ for real $x$, containing trigonometric and incomplete cylindrical functions.

\section*{Acknowledgments}
We are grateful to Prof. O.A. Chichigina for helpful and motivating discussions.

\bibliographystyle{amsplain}
\bibliography{asinc}

\end{document}